\documentclass[11pt]{amsart}
\usepackage{amssymb, latexsym}
\theoremstyle{plain}
\newtheorem{theorem}{Theorem}
\newtheorem{corollary}{Corollary}
\newtheorem*{thm-cheb}{Theorem (Chebyshev)}

\newtheorem{proposition}{Proposition}

\newtheorem*{2'}{Theorem 2'}
\newtheorem*{3'}{Theorem 3'}

\theoremstyle{remark}

\newtheorem*{Remark 1}{Remark 1}
\newtheorem*{Remark 2}{Remark 2}
\newtheorem*{Remark 3}{Remark 3}
\newtheorem*{Remark 4}{Remark 4}

\numberwithin{equation}{section}

\begin{document}

\title[  Secretary Problem with multiple items at each rank]
 {Two measures of efficiency for  the Secretary Problem with  multiple items at each rank }

\author{Ross G. Pinsky}


\address{Department of Mathematics\\
Technion---Israel Institute of Technology\\
Haifa, 32000\\ Israel}
\email{ pinsky@math.technion.ac.il}

\urladdr{https://pinsky.net.technion.ac.il/}

\subjclass[2000]{60G40, 60C05} \keywords{secretary problem, optimal stopping,  permutations}
\date{}

\begin{abstract}
For $2\le k\in\mathbb{N}$, consider the following adaptation of the classical secretary problem. There are $k$ items at each of $n$ linearly ordered ranks.
The $kn$ items are revealed, one item at a time, in a uniformly random order,
to an observer whose objective is to select an item of highest rank.
At each stage the observer only knows the relative ranks of the items that have arrived thus far, and must either select the current item, in which case the process terminates, or reject it and continue to the next item.
For $M\in\{0,1,\cdots, kn-1\}$, let $\mathcal{S}(n,k;M)$ denote the strategy whereby one  allows the first $M$ items to pass, and then selects the first later arriving item whose rank is \it either equal to or greater than\rm\ the highest rank of the first $M$ items (if such an item exists). Let $W_{\mathcal{S}(n,k;M)}$ denote the event that one selects an item of highest rank using  strategy $\mathcal{S}(n,k;M)$ and let
 $P_{n,k}(W_{\mathcal{S}(n,k;M)})$ denote the corresponding probability.
We obtain a formula for $P_{n,k}(W_{\mathcal{S}(n,k;M)})$, and  for $\lim_{n\to\infty}P_{n,k}(W_{\mathcal{S}(n,k;M_n)})$, when
$M_n\sim ckn$, with $c\in(0,1)$.
In the classical secretary problem ($k=1$), the asymptotically optimal strategy  $M_n\sim cn$ occurs with $c=\frac1e\approx0.368$, and the corresponding asymptotic probability of success is $\frac1e\approx 0.368$.
For $k=2$, the asymptotically optimal strategy  $M_n\sim ckn$ occurs with $c\approx 0.386$, and the limiting probability of success is about
0.701.
For $k=3$, the optimal probability
   is above 0.85, for $k=7$, that probability exceeds 0.99,  and for $k\ge12$, it is 1.000 to three decimal places.
In the problem with multiple items at each rank,  there is an  additional
measure of efficiency of a strategy besides the probability of selecting an item of highest rank; namely how quickly one selects an item of highest rank.
We give a rather complete picture of this efficiency.
We also consider the strategies $\mathcal{S}^+(n,k;M)$, for $M\in\{0,\cdots, kn-1\}$, whereby one allows the first $M$
items to pass, and then selects the first later arriving item whose rank is \it strictly greater than\rm\ the highest rank of the first $M$ items.
We show that
these   strategies  turn the problem   into one that is essentially equivalent to the classical
 secretary problem.

\end{abstract}

\maketitle
\section{Introduction and Statement of Results}
Recall the classical secretary problem: For $n\in\mathbb{N}$, a set of $n$ linearly ranked items is revealed, one item at a time, to an observer whose objective is to select the item
of highest rank. The order of the items is uniformly random; that is, each of the  $n!$ permutations of the ranks is equally likely.
At each stage, the observer only knows the relative ranks of the items that have arrived thus far, and must  either select the current item, in which case the process terminates, or reject it and continue to the next item. If the observer rejects the first $n-1$ items, then the $n$th and final item to arrive must be accepted. Since only relative ranks are known, the only reasonable strategies are $\{\mathcal{S}(n;M)\}_{M=0}^{n-1}$, where the strategy $\mathcal{S}(n;M)$ is to let the first $M$ items pass, and then to select the first
later arriving item that is ranked higher than all of the first $M$ items (if such an item exists).
As is very well known, asymptotically as $n\to\infty$, the optimal  strategy is $\mathcal{S}(n;M_n)$, where $M_n\sim \frac ne$.
The limiting probability of successfully selecting the item of highest rank is $\frac1e\approx0.368$.

In this paper, we extend the secretary problem to the case that there are multiple items at each rank. Fix an integer $2\le k\in\mathbb{N}$.
For  $n\in\mathbb{N}$, consider  $n$ linearly ordered ranks and  $kn$ items, with $k$ items at each rank.
The $kn$ items are revealed, one item at a time, to an observer whose objective is to select an item of highest rank. The order of the items is uniformly random.
At each stage the observer only knows the relative ranks of the items that have arrived thus far, and must either select the current item, in which case the process terminates, or reject it and continue to the next item.
Thus, the problem is equivalent to the problem of revealing the items one by one in  a uniformly random permutation of
the set $\cup_{i=1}^n\cup_{l=1}^k \{i_l\}$,
which consists of $k$ repetitions of each number $i\in[n]$.
The permutations of this set will be denoted by $S_{n,k}$. (Of course, $S_{n,k}$ is equivalent to $S_{nk}$, the set of permutations of $[nk]$.)
Each permutation in $S_{n,k}$ has $nk$ positions.
Here and in the sequel, the items are always revealed according to the left to right order of their positions in
the permutation. The items $\{n_l\}_{l=1}^k$ are considered the items of highest rank.

In the classical secretary problem, there is only one measure of efficiency; namely the probability of selecting the item of highest rank. In the problem at hand, there is an additional measure of efficiency; namely how quickly one selects an item of highest rank.  Indeed, as we shall see, this latter measure of efficiency becomes the dominant one
 because for all but the first few values of $k$, the main class of strategies we consider is capable  of  selecting an item of highest rank with probability extremely close to one.

The main results of this paper concern  a class of strategies of a form similar to, but slightly different from,  the above mentioned strategies in the classical secretary problem. For $0\le M<nk$,
 denote by $\mathcal{S}(n,k;M)$ the strategy whereby one  allows the first $M$ items to pass, and then selects the first later arriving item whose rank is \it either equal to or greater than\rm\ the highest rank of the first $M$ items (if such an item exists).

The paper \cite{GKM} considered a different and more intricate class  of strategies  for  the secretary problem considered here,
in the case that $k=2$, and  these strategies can in fact be defined for any $k$, as we now describe. For $0\le R\le n$, let $\mathcal{T}(n,k;R)$ denote the strategy whereby
one lets items pass until at least one item has occurred from $R$ of the $n$ different ranks. After that, one selects the first item to occur that satisfies the following two criteria:
(1) $k-1$ items of the same rank as the current item have already passed; and (2) the rank of the current item is larger than or equal to the rank of every item that
has preceded it.
If these two criteria never occur, then one  fails.

It seems intuitively clear that an optimal choice of strategies from the class $\mathcal{T}(n,k;R)$ will yield a higher probability of selecting an item of highest rank
than will an optimal choice of strategies from the class $\mathcal{S}(n,k;M)$. However, as already noted, for all but the first few values of $k$, the optimal
strategy from the class $\mathcal{S}(n,k;M)$ already yields probabilities extremely close to one. We will show   that the class of strategies $\mathcal{S}(n,k;M)$  is
overwhelmingly more efficient than the class $\mathcal{T}(n,k;R)$ from the point of view of how quickly an item of highest rank is selected.

 We also consider the  strategies $\mathcal{S}^+(n,k;M)$
  whereby one  allows the first $M$ items to pass, and then selects the first later arriving item whose rank is \it strictly greater than\rm\ the highest rank of the first $M$ items (if such an item exists).
  We shall see that the class $\mathcal{S}^+(n,k;M)$ of strategies  turns the problem with multiple items at each rank  into one that is essentially equivalent to the classical secretary problem.

We now turn to a presentation of the results.
  Let $P_{n,k}$ denote the uniform probability measure on the set of permutations   $S_{n,k}$. Let $W_{\mathcal{S}(n,k;M)}\subset S_{n,k}$
  denote the event that  an item of highest rank is selected when using the strategy $\mathcal{S}(n,k;M)$.
 Our first result gives an exact formula for  $P_{n,k}(W_{\mathcal{S}(n,k;M)})$.

 We use the notation $(b)_a=b(b-1)\cdots (b-a+1)$ for falling factorials, where
  $a,b\in\mathbb{Z}^+$.
 \begin{proposition}\label{prop1}
 \begin{equation}\label{formula1}
\begin{aligned}
& P_{n,k}(W_{\mathcal{S}(n,k;M)})=\sum_{l=1}^{k-1}\frac{\binom Ml(k)_l(k(n-1))_{M-l}}{(kn)_M}+\\
&k\sum_{j=1}^{n-1}\sum_{l=1}^k \frac{\binom Ml(k)_l(k(j-1))_{M-l}}{(kn)_M}\frac 1{k(n-j+1)-l}.
\end{aligned}
 \end{equation}
 \end{proposition}
  An asymptotic analysis of the formula in Proposition \ref{prop1} leads to the central result of this paper.
  \begin{theorem}\label{thm1}
  Let $M_n\sim ckn$, where $c\in(0,1)$. Then
  \begin{equation}\label{formula2}
  \begin{aligned}
  &\lim_{n\to\infty} P_{n,k}(W_{\mathcal{S}(n,k;M_n)})=-(1-c)^k\sum_{l=1}^{k-1}\binom kl(\frac c{1-c})^l\frac l{k-l}+\\
  &k\sum_{l=1}^{k-1}\binom klc^l\int_0^{1-c}\frac{y^{k-l-1}}{1-y^k}dy-c^k\log(1-(1-c)^k).
  \end{aligned}
  \end{equation}
  \end{theorem}
Using partial fractions and trigonometric substitution,   we can calculate explicitly the integrals on the right hand side above  for the cases $k=2,3$. We obtain
   \begin{equation}\label{k23}
   \begin{aligned}
&\lim_{n\to\infty}P_{n,2}(W_{\mathcal{S}(n,2;M_n)})=-2c(1-c)+(2c-c^2)\log(2-c)-(2c+c^2)\log c;\\
 &\lim_{n\to\infty}P_{n,3}(W_{\mathcal{S}(n,3;M_n)})=-\frac32(1-c)c(1+3c)-(3c+3c^2+c^3)\log c+\\
 &(\frac32c+\frac32c^2-c^3)\log(c^2-3c+3)+3\sqrt3(-c+c^2)\arctan(\frac{3-2c}{\sqrt3})+\frac{\sqrt3\pi}2(c-c^2).
   \end{aligned}
   \end{equation}
   \medskip

   Table 1  gives  the approximate optimal value of $c$ and the corresponding optimal limiting probability of selecting an item of highest rank for $k$ between 1 and 11.
   For $k=2$, the probability of success is about .701,  compared to $\frac1e\approx .368$  in the classical case when there is only one item of each rank, while
   the optimal choice of $c$, namely, $c\approx0.386$,  is close to the optimal value $\frac1e$ in the classical case.  For $k=3$, the optimal probability
   is above 0.85, for $k=7$, that probability exceeds 0.99,  and for $k=10$, it is approximately 0.999.

For $k\ge12$, the optimal probability of success is 1.000, when rounded off to three decimal places.
Table 2 considers several such values of $k$ and gives the approximate range of values of $c$ for which the probability of  selecting an item of highest rank is approximately equal to 1.000.

 \begin{table}[h!]
\centering
\begin{tabular}{||c c c ||}
 \hline
 k & argmax  for $c$& max. prob. \\ [0.5ex]
 \hline\hline
 2 & 0.386 & 0.701  \\
 3 & 0.413 & 0.854  \\
 4 & 0.431 & 0.928  \\
 5 & 0.444 & 0.964\\
 6 & 0.453 & 0.982 \\
 7 & 0.460 & 0.991  \\
 8 & 0.465 & 0.996 \\
 9 & 0.465 & 0.996\\
 10 & 0.472 & 0.999\\
 11&  0.474 & 0.999\\
 \hline
\end{tabular}
\\ [1ex]
\caption{Optimal $c$ and optimal probability}
\label{table:1}
\end{table}

 \begin{table}[h!]
\centering
\begin{tabular}{||c c c ||}
 \hline
 k & range for  $c$& probability \\ [0.5ex]
 \hline\hline
 12 & [.44,.52] & 1.000 \\
 15 & [.36,.60]& 1.000\\
 25 & [.24,.73] & 1.000 \\
 50 & [.13,.85] & 1.000  \\
 \hline
\end{tabular}
\\ [1ex]
\caption{Range of $c$ for which  probability $\approx1.000$}
\label{table:2}
\end{table}

As noted above, the paper \cite{GKM} considered the class of strategies $\mathcal{T}(n,k;R)$ in the case $k=2$.
 Let $W_{\mathcal{T}(n,k;R)}\subset S_{n,k}$
  denote the event that  an item of highest rank is selected when using the strategy $\mathcal{T}(n,k;R)$.
It was shown there that the limiting probability of success,
$\lim_{n\to\infty}P_{n,k}(W_{\mathcal{T}(n,2;R_n)})$, with $R_n\sim cn$, $c\in(0,1)$, is maximized for $c\approx0.4709$, and the corresponding limiting probability is approximately
$0.7680$. This probability is higher than the probability 0.701  obtained using strategy $\mathcal{S}(n,2;M_n)$ with $M_n\sim0.386(2n)$ as in Table 1.
It  isn't obvious  how to extend the analysis in \cite{GKM} to $k\ge3$.
On the other hand, in light of the results in Tables 1 and 2, for all but the first few values of $k$, a strategy from the class $\mathcal{T}(n,k;R)$  can hardly
be more effective at selecting an item of highest rank than is the optimal strategy from the class $\mathcal{S}(n,k;M)$.
Furthermore, we now show that the class of strategies $\mathcal{S}(n,k;M)$ has an overwhelming  advantage over the class of strategies
$\mathcal{T}(n,k;R)$. Indeed,  note that if strategy $\mathcal{T}(n,k;R)$ succeeds, then it automatically
selects the $k$th and final occurrence of an item with rank $n$, whereas strategy $\mathcal{S}(n,k;M_n)$ can succeed on an earlier occurrence of an item with rank $n$.
Taking this point of view, we now investigate and compare the  efficiency of these two classes of strategies.


On $W_{\mathcal{S}(n,k;M_n)}\subset S_{n,k}$,  the event  that an item of highest rank is selected when using strategy $\mathcal{S}(n,k;M_n)$, define $J_{n,k}^{(\mathcal{S})}$ with values in $[k]$ by $J_{n,k}^{(\mathcal{S})}=i$, if  the item of rank $n$ that was selected from among the $k$ items $\{n_l\}_{l=1}^k$ of rank $n$, is the $i$th item of rank $n$ to occur.
Also, on the event  $W_{\mathcal{S}(n,k;M_n)}$, let $X^{(\mathcal{S})}_{n,k}$
denote the
  position in $[kn]$ in which the selected item of highest rank occurs.
 (The dependence on $\{M_n\}$ in $J_{n,k}^{(\mathcal{S})}$ and $X^{(\mathcal{S})}_{n,k}$ is implicit in $\mathcal{S}=\mathcal{S}(n,k;M_n))$.

On  $W_{\mathcal{T}(n,k;R_n)}\subset S(n,k)$, the event  that an item of highest rank is selected when using strategy $\mathcal{T}(n,k;R_n)$,
 define $J_{n,k}^{(\mathcal{T})}=k$, because  from the definition of strategy
$\mathcal{T}$,
the item of highest rank $n$ that was selected is automatically
the $k$th and final  item of rank $n$ to occur.
Also, on the event $W_{\mathcal{T}(n,k;R_n)}$,
 let $X^{(\mathcal{T})}_{n,k}\in[kn]$ denote the
  position in which the selected item of highest rank occurs. (The dependence  on  $\{T_n\}$  in $X^{(\mathcal{T})}_{n,k}$ is implicit in $\mathcal{T}=\mathcal{T}(n,k;R_n)$.)

\begin{theorem}\label{1tok}
Using the strategy $\mathcal{S}(n,k;M_n)$,
where $M_n\sim ckn$, define
$$
p^{(k)}_c(i)=   \lim_{n\to\infty}P_{n,k}(J_{n,k}^{\mathcal{S}}=i|W_{\mathcal{S}(n,k;M_n)}),\ i\in[k].
$$
That is, $p_c^{(k)}(i)$ is the limiting probability as $n\to\infty$ that the selected    item of rank $n$  is the $i$th item of rank $n$ to occur.
Then
\begin{equation}\label{pi}
p_c ^{(k)}(i)=\begin{cases}\frac1{\lim_{n\to\infty} P_{n,k}(W_{\mathcal{S}(n,k;M_n)})}(1-c)^{k-i+1}c^{i-1}\binom k{i-1},\ i\in\{2,\cdots, k\};\\
1-\frac1{\lim_{n\to\infty} P_{n,k}(W_{\mathcal{S}(n,k;M_n)})}\big(1-c^k-(1-c)^k\big), \ i=1.
\end{cases}
\end{equation}
\end{theorem}
\bf\noindent Remark.\rm\    In particular, if $\lim_{n\to\infty} P_{n,k}(W_{\mathcal{S}(n,k;M_n)})\approx1$ with $M_n\sim c(kn)$ (see Table 1), then
$$
p_c^{(k)}(i)\approx\begin{cases}(1-c)^{k-i+1}c^{i-1}\binom k{i-1},\ i\in\{2,\cdots, k\};\\
c^k+(1-c)^k, \ i=1.
\end{cases}
$$
Equivalently, the distribution $\{p_c^{(k)}(i)\}_{i=1}^k$ is approximately the distribution of
the random variable $1+X_{k,c}1_{X_{k,c}\neq k}$, where $X_{k,c}$ has the binomial distribution with parameters $k$ and $c$.

Since $P_{n,k}(J_{n,k}^{\mathcal{T}}=k|W_{\mathcal{T}(n,k;R_n)})=1$, it follows that
in contrast to the above  distribution
$\{p_c ^{(k)}(i)\}_{i=1}^k$ on $[k]$ corresponding to the class of strategies $\mathcal{S}(n,k;M)$,  the analogous  distribution on $[k]$ corresponding to the class of strategies
$\mathcal{T}(n,k;R)$ is the degenerate distribution that places all its  probability on  $k\in[k]$.

The following theorem analyzes the quantity $\lim_{n\to\infty}\frac{E_{n,p}(X_{n,k}^{\mathcal{S}}|W_{\mathcal{S}(n,k;M_n)})}{kn}$, the
limiting proportion of items that are observed until an item of highest rank $n$ is selected, conditioned on an item of highest rank begin selected,
when using the strategy $\mathcal{S}(n,k;M_n)$.
\begin{theorem}\label{1tokn}
Using the strategy $\mathcal{S}(n,k;M_n)$, where $M_n\sim ckn$, define
$$
\gamma^{(k)}(c)=\lim_{n\to\infty}\frac{E_{n,p}(X_{n,k}^{\mathcal{S}}|W_{\mathcal{S}(n,k;M_n)})}{kn}.
$$
That is, $\gamma^{(k)}(c)$ is the limiting proportion of items that are observed until an item of rank $n$ is selected.
Then
\begin{equation}\label{gamma}
\begin{aligned}
&\gamma^{(k)}(c)=c+(1-c)\sum_{i=1}^k \frac{p_c^{(k)}(i)}{k+2-i}=\\
&c+\frac{1-c}{\lim_{n\to\infty}P_{n,k}(W_{\mathcal{S}(n,k;M_n)})}E\frac1{k+1-X_{k,c}}1_{X_{k,c}\not\in \{0,k\}}+\\
&\frac{1-c}{k+1}\big(1-\frac1{\lim_{n\to\infty}P_{n,k}(W_{\mathcal{S}(n,k;M_n)})}(1-c^k-(1-c)^k)\big),
\end{aligned}
\end{equation}
where $p_c^{(k)}(l)$ is as in \eqref{pi} and $X_{k,c}$ is a binomial random variable with parameters $k$ and $c$.
\end{theorem}
\bf\noindent Remark.\rm\ In cases where $\lim_{n\to\infty}P_{n,k}(W_{\mathcal{S}(n,k;M_n)})\approx1$,
\eqref{gamma} reduces to
\begin{equation}\label{rem}
\begin{aligned}
&\gamma^{(k)}(c)\approx
c+(1-c)\sum_{i=1}^k(1-c)^{k-i+1}c^{i-1}\binom k{i-1}\frac1{k+2-i}+\frac{(1-c)c^k}{k+1}=\\
& c+(1-c)E\frac1{k+1-X_{k,c}1_{X_{k,c}\neq k}}.
\end{aligned}
\end{equation}
By an elementary  large deviations result, $P(X_{k,c}\le (c-\epsilon)k)$ is exponentially small in $k$, for any $\epsilon>0$. From this, we obtain immediately the following corollary of \eqref{gamma}.
\begin{corollary}
For any sequence $\{c_k\}_{k=2}^\infty$ for which
$$
\liminf_{k\to\infty}\lim_{n\to\infty}P_{n,k}(W_{\mathcal{S}(n,k;M_n)})|_{M_n\sim c_kkn})>0
$$
and $\limsup_{k\to\infty}c_k<1$, one has
$$
\gamma^{(k)}(c_k)=c_k+\theta(\frac1k),\ \text{as}\ k\to\infty.
$$
\end{corollary}


 \begin{table}[h!]
\centering
\begin{tabular}{||c c c ||}
 \hline
 k & argmax  for $c$ from Table 1& $\gamma^{(k)}(c)$ \\ [0.5ex]
 \hline\hline
 2 & 0.386 & 0.660  \\
 3 & 0.413 & 0.636  \\
 4 & 0.431 & 0.618  \\
 5 & 0.444 & 0.604\\
 6 & 0.453 & 0.592 \\
 7 & 0.460 & 0.583  \\
 8 & 0.465 & 0.575 \\
 9 & 0.465 & 0.565\\
 10 & 0.472 & 0.563\\
 \hline
\end{tabular}
\\ [1ex]
\caption{Value of $\gamma^{(k)}(c)$, the limiting proportion of items that are observed until an item of rank $n$ is selected,
 at optimal $c$ from Table 1}
\label{table:3}
\end{table}

 \begin{table}[h!]
\centering
\begin{tabular}{||c c c ||}
 \hline
 k & lower end approx. argmax for $c$& $\gamma^{(k)}(c)$ \\ [0.5ex]
 \hline\hline
 12 & .44 & .52 \\
 15 & .36& .42\\
 25 & .24 &  .28\\
 50 & .13 &   .15\\
 \hline
\end{tabular}
\\ [1ex]
\caption{Value of $\gamma^{(k)}(c)$, the limiting proportion of items that are observed until an item of rank $n$ is selected, for the smallest value of $c$  for which  probability $\approx1.000$}
\label{table:4}
\end{table}

Table 3 gives the value of $\gamma^{k)}(c)$, the limiting
proportion of items that are observed until an item of highest rank is selected, conditioned on an item of highest rank being selected, for $2\le k\le 10$, with the optimal value of $c$ from Table 1.
Table  4 considers several values of $k$ with $k\ge12$, the range of $k$ for which the optimal probability is 1.000, to three decimal places. Numerical results show that
the function in \eqref{rem} is increasing in $c$. Thus,
values are given for $\gamma^{(k)}(c)$
using  the lowest value of $c$, taken from Table 2, for which the probability of success is approximately  1.000.

With regard to the limiting proportion of  items that are observed until an item of rank $n$ is selected,
in contrast to Theorem \ref{1tokn} and Tables 3 and 4 for the class of strategies $\mathcal{S}(n,k;M)$, we have the following result for the class of strategies
$\mathcal{T}(n,k;R)$.
\begin{proposition}\label{1toknT}
Using any sequence of strategies $\{\mathcal{T}(n,k;R_n)\}$ for which $\alpha_k:=\alpha_k(\{R_n\})=\lim_{n\to\infty}P_{n,k}(W_{\mathcal{T}(n,k;R_n)})$ exists, one has
\begin{equation}\label{knT}
\begin{aligned}
&\frac1{\alpha_k}\frac k{k+1}-\frac{1-\alpha_k}{\alpha_k}
\le
\liminf_{n\to\infty}E_{n,k}(X_{n,k}^{\mathcal{T}}|W_{\mathcal{T}(n,k;R_n)})\le\\
&\limsup_{n\to\infty}E_{n,k}(X_{n,k}^{\mathcal{T}}|W_{\mathcal{T}(n,k;R_n)})\le
\frac1{\alpha_k}\frac k{k+1}.
\end{aligned}
\end{equation}
\end{proposition}
\noindent\bf Remark.\rm\ In particular, it follows that if $\lim_{n\to\infty}P_{n,k}(W_{\mathcal{T}(n,k;R_n)})\approx1$, then
$\lim_{n\to\infty}E_{n,k}(X_{n,k}^{\mathcal{T}}|W_{\mathcal{T}(n,k;R_n)})\approx \frac k{k+1}$.


\medskip

Our final results concern  the class of strategies  $\mathcal{S}^+(n,k;M)$ that was defined above after the definitions of the classes of   strategies $\mathcal{S}(n,k;M)$ and
$\mathcal{T}(n,k;R)$. As with the strategies $\mathcal{S}(n,k;M)$, we give an exact result for fixed $n$ and an asymptotic result as $n\to\infty$.
 Let $W_{\mathcal{S}^+(n,k;M)}\subset S_{n,k}$
  denote the event that  an item of highest rank is selected when using the strategy $\mathcal{S}^+(n,k;M)$.
\begin{proposition}\label{prop+}
\begin{equation}\label{formula1+}
\begin{aligned}
 P_{n,k}(W_{\mathcal{S}^+(n,k;M)})=\sum_{j=1}^{n-1}\sum_{l=1}^k \frac{\binom Ml(k)_l(k(j-1))_{M-l}}{(kn)_M}\frac 1{n-j}.
\end{aligned}
 \end{equation}
\end{proposition}
\begin{theorem}\label{thm+}
 Let $M_n\sim ckn$, where $c\in(0,1)$. Then
  \begin{equation}\label{formula2+}
  \begin{aligned}
  &\lim_{n\to\infty} P_{n,k}(W_{\mathcal{S}^+(n,k;M_n)})=-\big(1-(1-c)^k\big)\log(1-(1-c)^k).
  \end{aligned}
  \end{equation}
  \end{theorem}
As is very well known, for the classical secretary problem with strategy $\mathcal{S}(n;M)$, as defined in the first paragraph of the paper, the asymptotic limiting probability of selecting the item of highest rank is equal to
$-c\log c$, if $M_n\sim cn$, $c\in(0,1)$. Thus, Theorem \ref{thm+} shows that the asymptotic limiting probability
of selecting an item of highest rank in  the secretary  problem with multiple items at each rank when using the strategy $\mathcal{S}^+(n,k;M_n)$
with $M_n\sim c^+kn$ is equal to the asymptotic limiting probability of selecting the item of highest rank in the classical secretary problem when using the strategy $\mathcal{S}(n;M_n)$
with $M_n\sim cn$, where $c^+=1-(1-c)^\frac1k$. Since $-c\log c$, for $c\in(0,1)$, attains its maximum value of $\frac1e$ at $x=\frac1e$, the following corollary of Theorem \ref{thm+} is immediate.
\begin{corollary}
The asymptotically optimal strategy from among the strategies $\mathcal{S}^+(n,k;M_n)$ is the one with
$M_n\sim \big(1-(1-\frac1e)^\frac1k\big)kn$, and the corresponding optimal limiting probability of selecting an item of highest rank is $\frac1e$.
\end{corollary}
The reason   the secretary problem with multiple items  at each rank, considered with the
strategies $\mathcal{S}^+(n,k;M)$,  is essentially turned into a classical secretary problem
is that under these strategies the only highest ranked item that  is possible to select
is the first one that occurs.
\medskip

The classical secretary problem, where there is one item at each rank, but  with adaptations to increase the chance of winning, go all the way back to the fundamental  paper of Gilbert and Mosteller \cite{GM}.
In particular, they considered the situation where one is given $r$ opportunities  to select the highest ranked item, as well as the situation in which one is given one opportunity to  select
an item from among the $r$ top ranked items.
In the first situation above, when there are  $r$ opportunities to select the highest ranked item, they showed that the asymptotic probability of winning when using the best strategy
is about $0.591$ for $r=2$, and increases to about $0.965$ when $r=8$.  In the second situation above, when there is one opportunity to select an item from among the $r$ top ranked items, they showed that the  asymptotic probability of winning when using the best strategy is about 0.574 for $r=2$. The authors did not analyze  their formula numerically  in  cases with $r>2$.

For the  secretary problem  in its classical setup, but  with items arriving in a non-uniform order, see for example \cite{GM, HK, Pf} and the recent paper \cite{P1}.
See \cite{B} for another approach to the secretary problem, and see \cite{F} for a review of a variety of other variations on the secretary problem.

We prove
Proposition \ref{prop1} in section \ref{sec2} and Theorem \ref{thm1} in section \ref{sec3}.
In section \ref{aux}, we prove a couple of auxiliary results that are needed for the proofs of Proposition \ref{1toknT} and Theorem \ref{1tokn}, and then give the proof of
Proposition \ref{1toknT}.
In section \ref{quick1} we prove Theorem \ref{1tok} and in section \ref{quick2} we prove Theorem \ref{1tokn}. The proofs of these two theorems rely heavily on the calculations in the proofs of Proposition \ref{prop1} and Theorem \ref{thm1}.
The proofs of Proposition \ref{prop+} and Theorem \ref{thm+} are obtained quickly in section
\ref{sec4} by making small changes in the proofs of Proposition \ref{prop1} and Theorem \ref{thm1}.

\section{Proof of Proposition \ref{prop1}}\label{sec2}
Fix $2\le k\in\mathbb{N}$. For $n\in\mathbb{N}$, we consider a uniform permutation of the multi-set $\{1^k,2^k,\cdots, n^k\}$. The $kn$ items of this multi-set arrive according to the order of this permutation.
We consider $n$ to be the highest rank and 1 to be the lowest rank.
We wish to calculate the probability $P_{n,k}(\mathcal{S}(n,k;M))$ of winning under the strategy $\mathcal{S}(n,k;M)$, which was introduced in the third paragraph of the paper.

Let $A^{(n)}_{M,j,l}$ denote the event that among the first $M$ items, the number $j$ has occurred $l$ times, and no number greater than $j$ has occurred at all,
 where $1\le M\le kn-1,\ 1\le j\le n$ and  $l\in[k]$.
To calculate $P_{n,k}(A^{(n)}_{M,j,l})$, it will be convenient to consider all  $nk$ objects as distinguishable from one another. (For each $i\in[n]$, think of the $k$ different $i$'s as $\{i_l\}_{l=1}^k$.)
Then there are $(kn)_M$ different possible outcomes for the first $M$ items. In order for $A^{(n)}_{M,j,l}$ to occur, one needs to select $l$ locations to place the  number $j$, one needs to select $l$ $j$'s from among the
$k$ different $j$'s, and one needs to
fill the other $M-l$ locations with numbers less than $j$.  Thus,
\begin{equation}\label{A}
P_{n,k}(A^{(n)}_{M,j,l})=
\frac{\binom Ml(k)_l(k(j-1))_{M-l}}{(kn)_M}.
\end{equation}
(Of course, the above probability is zero if $M-l>k(j-1)$.)

Consider now the case that $j\in[n-1]$. If one employs the strategy $\mathcal{S}(n,k;M)$, and the event $A^{(n)}_{M,j,l}$ occurs, then after the $M$th items arrives, one will select the first item larger or equal to $j$. Among the items arriving after
the $M$th item arrives, there are $k(n-j+1)-l$ items that are larger or equal to $j$, of which $k$ of them have the highest rank $n$. Since items arrive in uniformly random order,
\begin{equation}\label{condprobj<n}
P_{n,k}(W_{\mathcal{S}(n,k;M)}|A^{(n)}_{M,j,l})=\frac k{k(n-j+1)-l},\ j\in[n-1].
\end{equation}

 On the other hand, consider now the case  $j=n$. If one  employs the strategy $\mathcal{S}(n,k;M)$, and the event $A^{(n)}_{M,n,l}$ occurs, then after the $M$th item arrives,  one will
 select the first item that is equal to  $n$, if such an item exists. Of course, it will exist if $l\in[k-1]$ and it won't exist if $l=k$.
 Thus,
 \begin{equation}\label{condprob=n}
 P_{n,k}(W_{\mathcal{S}(n,k;M)}|A^{(n)}_{M,n,l})=\begin{cases} 1,\ \text{if}\ l\in[k-1];\\ 0,\  \text{if}\ l=k.\end{cases}
 \end{equation}

 Since for each $M$,  the collection of events $\{A^{(n)}_{M,j,l}: j\in[n], l\in [k]\thinspace\}$ are disjoint, and the probability of their union is 1, it follows from \eqref{A}-\eqref{condprob=n} that
 \begin{equation*}
\begin{aligned}
& P_{n,k}(W_{\mathcal{S}(n,k;M)})=\sum_{l=1}^{k-1}\frac{\binom Ml(k)_l(k(n-1))_{M-l}}{(kn)_M}+\\
&k\sum_{j=1}^{n-1}\sum_{l=1}^k \frac{\binom Ml(k)_l(k(j-1))_{M-l}}{(kn)_M}\frac 1{k(n-j+1)-l}.
\end{aligned}
 \end{equation*}
 \hfill $\square$

\section{Proof of Theorem \ref{thm1}}\label{sec3}
Let  $M_n=c_nkn$, be an integer for each $n\in\mathbb{N}$, with $\lim_{n\to\infty}c_n=c$. For $N\in\mathbb{N}$,
we write the second term on the right hand side of  \eqref{formula1}, with $M=M_n$, as
\begin{equation}\label{secondterm}
\begin{aligned}
&k\sum_{j=1}^{n-1}\sum_{l=1}^k \frac{\binom{M_n}l(k)_l(k(j-1))_{M_n-l}}{(kn)_{M_n}}\frac 1{k(n-j+1)-l}=\\
&k\sum_{l=1}^k\binom{c_nkn}l(k)_l\sum_{j=1}^{n-N}\frac{(k(j-1))_{c_nkn-l}}{(kn)_{c_nkn}}\frac 1{k(n-j+1)-l}+\\
&k\sum_{l=1}^k\binom{c_nkn}l(k)_l\sum_{j=n-N+1}^{n-1}\frac{(k(j-1))_{c_nkn-l}}{(kn)_{c_nkn}}\frac 1{k(n-j+1)-l}.
\end{aligned}
\end{equation}
For sufficiently large $n$, there exists a constant $C=C(k,l,c)$ such that
\begin{equation}\label{est1N}
\binom{c_nkn}l\frac{(k(j-1))_{c_nkn-l}}{(kn)_{c_nkn}}\le C\frac{(k(j-1))_{c_nkn-l}}{(kn)_{c_nkn-l}}\le C\frac{(kj)_{c_nkn-l}}{(kn)_{c_nkn-l}}.
\end{equation}
Using the fact that $\frac{a-x}{b-x}$ is decreasing in $x$ for $0<x<a<b$, it follows that for any $c'\in(0,c)$ and for sufficiently large $n$,
\begin{equation}\label{est2N}
\frac{(kj)_{c_nkn-l}}{(kn)_{c_nkn-l}}\le (\frac{kj}{kn})^{c_nkn-l}\le (\frac jn)^{c'kn},\ l\in[k].
\end{equation}
From \eqref{est1N} and \eqref{est2N} we have for large $n$,
\begin{equation}\label{estN3}
\binom{c_nkn}l\sum_{j=1}^{n-N}\frac{(k(j-1))_{c_nkn-l}}{(kn)_{c_nkn}}\frac 1{k(n-j+1)-l}\le C\sum_{j=1}^{n-N}(\frac jn)^{c'kn}.
\end{equation}
Also,
\begin{equation}\label{estN4}
\begin{aligned}
&\sum_{j=1}^{n-N}(\frac jn)^{c'kn}=n\sum_{j=1}^{n-N}\frac1n(\frac jn)^{c'kn}\le n\int_0^{1-\frac {N-1}n}x^{c'kn}dx\le\\
& \frac n{c'kn+1}(1-\frac{N-1}n)^{c'kn+1}\le \frac n{c'kn+1}e^{-\frac{N-1}n(c'kn+1)}.
\end{aligned}
\end{equation}
From \eqref{estN3} and \eqref{estN4}, we conclude that the expression on the second line of \eqref{secondterm} satisfies
\begin{equation}\label{Nnegligible}
\lim_{N\to\infty}\limsup_{n\to\infty}k\sum_{l=1}^k\binom{c_nkn}l(k)_l\sum_{j=1}^{n-N}\frac{(k(j-1))_{c_nkn-l}}{(kn)_{c_nkn}}\frac 1{k(n-j+1)-l}=0.
\end{equation}

Letting $s=n-j$, we write the expression on the third line of \eqref{secondterm} as
\begin{equation}\label{keyterm}
\begin{aligned}
&k\sum_{l=1}^k\binom{c_nkn}l(k)_l\sum_{j=n-N+1}^{n-1}\frac{(k(j-1))_{c_nkn-l}}{(kn)_{c_nkn}}\frac 1{k(n-j+1)-l}=\\
&k\sum_{l=1}^k\binom{c_nkn}l(k)_l\sum_{s=1}^{N-1}\frac{(k(n-s-1))_{c_nkn-l}}{(kn)_{c_nkn}}\frac 1{k(s+1)-l}.
\end{aligned}
\end{equation}
We write
\begin{equation}\label{productanalysis1}
\binom{c_nkn}l\frac{(k(n-s-1))_{c_nkn-l}}{(kn)_{c_nkn}}=\frac{\binom{c_nkn}l}{\prod_{i=1}^l(kn-c_nkn+i)}
\prod_{t=0}^{c_nkn-l-1}\frac{k(n-s-1)-t}{kn-t}.
\end{equation}
We have
\begin{equation}\label{productanalysis2}
\prod_{t=0}^{c_nkn-l-1}\frac{k(n-s-1)-t}{kn-t}=\prod_{i=(1-c_n)kn+l+1}^{kn} (1-\frac{k(s+1)}i),
\end{equation}
and
\begin{equation}\label{log}
\begin{aligned}
&\log\prod_{i=(1-c_n)kn+l+1}^{kn} (1-\frac{k(s+1)}i)=\sum_{i=(1-c_n)kn+l+1}^{kn} \log(1-\frac{k(s+1)}i)=\\
&-k(s+1)\sum_{i=(1-c_n)kn+l+1}^{kn}\frac1i+o(1)=k(s+1)\log(1-c_n)+o(1),\ \text{as}\ n\to\infty.
\end{aligned}
\end{equation}
From \eqref{productanalysis2} and \eqref{log}, we obtain
\begin{equation}\label{productanalysis3}
\prod_{t=0}^{c_nkn-l-1}\frac{k(n-s-1)-t}{kn-t}=(1+o(1))(1-c)^{k(s+1)},\ \text{as}\ n\to\infty.
\end{equation}
Also, we have
\begin{equation}\label{productanalysis4}
\lim_{n\to\infty}\frac{\binom{c_nkn}l}{\prod_{i=1}^l(kn-c_nkn+i)}=\frac1{l!}(\frac c{1-c})^l.
\end{equation}
Now \eqref{productanalysis1}, \eqref{productanalysis3} and \eqref{productanalysis4} yield
\begin{equation}\label{keylimit}
\lim_{n\to\infty}\binom{c_nkn}l\frac{(k(n-s-1))_{c_nkn-l}}{(kn)_{c_nkn}}=\frac1{l!}(\frac c{1-c})^l(1-c)^{k(s+1)}.
\end{equation}
Letting $n\to\infty$ in \eqref{keyterm}, it follows from \eqref{keylimit} that the expression on the third line of \eqref{secondterm} satisfies
\begin{equation}\label{keytermlimit}
\begin{aligned}
&\lim_{n\to\infty}k\sum_{l=1}^k\binom{c_nkn}l(k)_l\sum_{j=n-N+1}^{n-1}\frac{(k(j-1))_{c_nkn-l}}{(kn)_{c_nkn}}\frac 1{k(n-j+1)-l}=\\
&k\sum_{l=1}^k\binom kl(\frac c{1-c})^l\sum_{s=1}^{N-1}(1-c)^{k(s+1)}\frac 1{k(s+1)-l}.
\end{aligned}
\end{equation}
From \eqref{secondterm},\eqref{Nnegligible} and \eqref{keytermlimit}, we conclude that the second term on the right hand side of \eqref{formula1}, with $M=M_n$,  satisfies
\begin{equation}\label{finalsecondtermformula}
\begin{aligned}
&\lim_{n\to\infty}k\sum_{j=1}^{n-1}\sum_{l=1}^k \frac{\binom{M_n}l(k)_l(k(j-1))_{M_n-l}}{(kn)_{M_n}}\frac 1{k(n-j+1)-l}=\\
&k\sum_{l=1}^k\binom kl(\frac c{1-c})^l\sum_{s=1}^\infty(1-c)^{k(s+1)}\frac 1{k(s+1)-l},\ \text{for}\  M_n\sim cn, \ c\in(0,1).
\end{aligned}
\end{equation}

We now consider the first term on the right hand side of \eqref{formula1}, with $M=M_n=c_nkn\sim ckn$.
 Although we considered
\eqref{keylimit} for $s\ge1$, of course it also holds for $s=0$. Thus, the first term on the right hand side of \eqref{formula1} satisfies
\begin{equation}\label{finalfirsttermformula}
\lim_{n\to\infty}\sum_{l=1}^{k-1}\frac{\binom {M_n}l(k)_l(k(n-1))_{M_n-l}}{(kn)_{M_n}}=
(1-c)^k\sum_{l=1}^{k-1}\binom kl(\frac c{1-c})^l.
\end{equation}

From \eqref{formula1}, \eqref{finalsecondtermformula} and \eqref{finalfirsttermformula}, we obtain
\begin{equation}\label{answer1}
\begin{aligned}
&\lim_{n\to\infty} P_{n,k}(\mathcal{S}(n,k;M_n))=k\sum_{l=1}^k\binom kl(\frac c{1-c})^l\sum_{s=1}^\infty(1-c)^{k(s+1)}\frac 1{k(s+1)-l}+\\
&(1-c)^k\sum_{l=1}^{k-1}\binom kl(\frac c{1-c})^l,\ \text{for}\ M_n\sim cn,\ c\in(0,1).
\end{aligned}
\end{equation}
We now analyze the infinite series
on the right hand side of \eqref{answer1}.

Let
\begin{equation}\label{G}
G_{k,l}(x)=\sum_{s=1}^\infty\frac{x^{k(s+1)}}{k(s+1)-l}.
\end{equation}
Then
$$
(x^{-l}G_{k,l})'=\sum_{s=1}^\infty x^{k(s+1)-l-1}=\frac{x^{2k-l-1}}{1-x^k}=\begin{cases}-x^{k-l-1}+\frac{x^{k-l-1}}{1-x^k},\ l=1,\cdots, k-1;\\ \frac{x^{k-1}}{1-x^k},\ l=k.\end{cases}
$$
Integrating and noting that $G_{k,l}$ vanishes at zero to the order $2k$, we obtain
\begin{equation}\label{Gformula}
G_{k,l}(x)=\begin{cases} -\frac{x^k}{k-l}+x^l\int_0^x\frac{y^{k-l-1}}{1-y^k}dy,\ l=1,\cdots, k-1;\\-\frac{x^k}k\log(1-x^k),\ l=k.\end{cases}
\end{equation}
From \eqref{G} and \eqref{Gformula}, we can rewrite the first term on the right hand side of \eqref{answer1} as
\begin{equation}\label{termI}
\begin{aligned}
&k\sum_{l=1}^k\binom kl(\frac c{1-c})^l\sum_{s=1}^\infty(1-c)^{k(s+1)}\frac 1{(s+1)k-l}=\\
&-k(1-c)^k\sum_{l=1}^{k-1}\binom kl(\frac c{1-c})^l\frac1{k-l}+k\sum_{l=1}^{k-1}\binom klc^l\int_0^{1-c}\frac{y^{k-l-1}}{1-y^k}dy-\\
&c^k\log(1-(1-c)^k).
\end{aligned}
\end{equation}
Now \eqref{formula2} follows from \eqref{answer1} and \eqref{termI}.
\hfill $\square$

\section{Two auxiliary results and the proof of Proposition \ref{1toknT}.}\label{aux}
We state and  prove two auxiliary results. The first one will be used in the proof of the second one  and in the proof  of Proposition \ref{1toknT},  and the second one
will be used in the proof of Theorem \ref{1tokn}.  We end this section with the proof of  Proposition \ref{1toknT}.

On the probability space $(S_{n,k},P_{n,k})$, for each $i\in[k]$, let $Y^{(n,k)}_i$ denote the position in $[kn]$ in which occurs the $i$th item to occur from among the $k$ items $\{n_l\}_{l=1}^k$ of
rank $n$.
(For example, if $n=4$ and $k=2$, and the permutation in $S_{4,2}$ is given by
$2_14_23_23_11_12_24_11_2$, then on this permutation, $Y^{(4,2)}_1$ takes  the value 2 and  $Y^{(4,2)}_2$ takes the value 7.)
\begin{proposition}\label{combi}
For $k\in\mathcal{N}$ and $i\in[k]$,
\begin{equation}\label{combiSnk}
\lim_{n\to\infty}\frac{E_{n,k}Y^{(n,k)}_i}{nk}=\frac i{k+1}.
\end{equation}
\end{proposition}
\begin{proof}
Since we will only need the cases $i=1$ and $i=k$, and since the latter case follows from the former one by considerations of symmetry, we will only proof the case
$i=1$. The other cases can be proved similarly.
The random variable $\frac{Y^{(n,k)}_1}{nk}$ on the probability space $(S_{n,k},P_{n,k})$ is a discretization of the random variable
$\min(U_1,\cdots, U_k)$, where $\{U_j\}_{j=1}^k$ are IID uniform random variables on $[0,1]$. Thus
$$
\lim_{n\to\infty}\frac{E_{n,k}Y^{(n,k)}_1}{nk}=E\min(U_1,\cdots, U_k).
$$
Since $P(\min(U_1,\cdots, U_k)\ge x)=(1-x)^k$, the density of $\min(U_1,\cdots, U_k)$ is $k(1-x)^{k-1}$. Thus,
$$
\begin{aligned}
&E\min(U_1,\cdots, U_k)=\int_0^1x\big(k(1-x)^{k-1}\big)dx=k\int_0^1(1-x)x^{k-1}dx=\\
&k(\frac1k-\frac1{k+1})=\frac1{k+1}.
\end{aligned}
$$
\end{proof}
\begin{proposition}\label{combibox}
A box contains $N$ black balls and   $A$ red balls.
Draw the balls from  the box without replacement until a red ball is obtained, and let $\tau^{(N,A)}$ denote the number of draws required.
Then
\begin{equation}\label{boxballs}
\lim_{N\to\infty}\frac{E\tau^{(N,A)}}N=\frac1{A+1}.
\end{equation}
\end{proposition}
\begin{proof}
Consider drawing all of the balls from the box, where the black balls are numbered from 1 to $N$ and the red balls are numbered from 1 to $A$. Every particular
order in which the balls are drawn is a permutation of the $N+A$ items.
Thus, if $N$ is a multiple of $A$, say $N=MA$, for $M\in\mathbb{N}$,  then $\tau^{(N,A)}$ is equivalent to $Y^{(n,k)}_1$ in Proposition \ref{combiSnk}, with $n=M+1$ and $k=A$.
Therefore, from Proposition \ref{combiSnk}, $\lim_{M\to\infty}\frac{E\tau^{(MA,A)}}{(M+1)A}=\frac1{A+1}$, or equivalently,
\begin{equation}\label{MA}
\lim_{M\to\infty}\frac{E\tau^{(MA,A)}}{MA}=\frac1{A+1}.
\end{equation}
A simple argument extends \eqref{MA} to the case that $MA$ is replaced by a generic $N$.
\end{proof}

\noindent \it Proof of Proposition \ref{1toknT}.\rm\
As in the statement of the proposition, let $\alpha_k=\lim_{n\to\infty}P_{n,k}(W_{\mathcal{T}(n,k;R_n)})$.
From the definition of the strategy  $\mathcal{T}(n,k;R_n)$, it follows that
$X_{n,k}^{\mathcal{T}}=Y^{(n,k)}_k|_{W_{\mathcal{T}(n,k;R_n)}}$
  where $Y^{(n,k)}_k$ is as defined before Proposition \ref{combi}.
  (Recall that $X_{n,k}^{\mathcal{T}}$  is only defined on $W_{\mathcal{T}(n,k;R_n)}$.)
Thus,
\begin{equation}\label{upperbdn}
E_{n,k}(X_{n,k}^{\mathcal{T}}|W_{\mathcal{T}(n,k;R_n)})=E_{n,k}(Y^{(n,k)}_k|W_{\mathcal{T}(n,k;R_n)})\le
\frac{E_{n,k}Y^{(n,k)}_k}{P_{n,k}(W_{\mathcal{T}(n,k;R_n)})},
\end{equation}
and then from \eqref{upperbdn} and \eqref{combiSnk}
\begin{equation}\label{upperbd}
\limsup_{n\to\infty}\frac1{nk}E_{n,k}(X_{n,k}^{\mathcal{T}}|W_{\mathcal{T}(n,k;R_n)})\le \frac1{\alpha_k}\frac k{k+1}.
\end{equation}
Using the notation $E(Z;B):=EZ1_B$, where $Z$ is a random variable and $B$ is an event, we also have
\begin{equation}\label{lowerbdn}
\begin{aligned}
&E_{n,k}(X_{n,k}^{\mathcal{T}}|W_{\mathcal{T}(n,k;R_n)})=E_{n,k}(Y^{(n,k)}_k|W_{\mathcal{T}(n,k;R_n)})=
\frac{E_{n,k}(Y^{(n,k)}_k;W_{\mathcal{T}(n,k;R_n)})}{P_{n,k}(W_{\mathcal{T}(n,k;R_n)})}\ge\\
&\frac{E_{n,k}Y^{(n,k)}_k}{P_{n,k}(W_{\mathcal{T}(n,k;R_n)})}
-nk\frac{1-P_{n,k}(W_{\mathcal{T}(n,k;R_n)})}{P_{n,k}(W_{\mathcal{T}(n,k;R_n)})},
\end{aligned}
\end{equation}
and then from \eqref{lowerbdn} and \eqref{combiSnk},
\begin{equation}\label{lowerbd}
\liminf_{n\to\infty}\frac1{nk}E_{n,k}(X_{n,k}^{\mathcal{T}}|W_{\mathcal{T}(n,k;R_n)})\ge\frac1{\alpha_k}\frac k{k+1}-\frac{1-\alpha_k}{\alpha_k}.
\end{equation}
Now  \eqref{upperbd} and \eqref{lowerbd} give \eqref{knT}.
\hfill $\square$

\section{Proof of Theorem \ref{1tok}}\label{quick1}
\it\noindent  Proof of Theorem \ref{1tok}.\rm\ Recall from the proof of Proposition \ref{prop1} that
$A^{(n)}_{M,j,l}$ denotes the event that among the first $M$ items,
 the number $j$ has occurred $l$ times, and no number greater than $j$ has occurred at all,
 where $1\le M\le kn-1,\ j\in[n]$ and  $l\in[k]$.
If $A^{(n)}_{M,j,l}$ occurs, for some $l\in[k]$ and some  $j\in[n-1]$, and  the strategy $\mathcal{S}(n,k;M)$ succeeds, then by the definition of the strategy,
$J_{n,k}^\mathcal{S}=1$. By \eqref{condprobj<n}, the probability that the strategy $\mathcal{S}(n,k;M)$ succeeds, conditioned on such an $A^{(n)}_{M,j,l}$, is equal to
$\frac k{k(n-j+1)-l}$.

On the other hand, if $A^{(n)}_{M,n,l}$ occurs, for some $l\in[k-1]$,
then by
the definition of the strategy, $J_{n,k}^\mathcal{S}=l+1$. Conditioned on $A^{(n)}_{M,n,l}$, it follows by
 \eqref{condprob=n} that the strategy $\mathcal{S}(n,k;M)$ automatically succeeds.
 Conditioned on $A^{(n)}_{M,n,k}$, it follows  by \eqref{condprob=n} that  the strategy $\mathcal{S}(n,k;M)$ automatically fails.

From the considerations in the above two paragraphs, along with \eqref{A}, we have
\begin{equation}\label{lge2}
\begin{aligned}
&P_{n,k}(J_{n,k}^{\mathcal{S}}=l|W_{\mathcal{S}(n,k;M_n)})=P_{n,k}(A^{(n)}_{M_n,n,l-1}|W_{\mathcal{S}(n,k;M_n)})=
\frac{P_{n,k}(A^{(n)}_{M_n,n,l-1})}{P_{n,k}(W_{\mathcal{S}(n,k;M_n)})}=\\
&\frac1{P_{n,k}(W_{\mathcal{S}(n,k;M_n)})}\frac{\binom {M_n}{l-1}(k)_{l-1}(k(n-1))_{M_n-l+1}}{(kn)_{M_n}},\ l\in\{2,\cdots, k\}.
\end{aligned}
\end{equation}
Writing $M_n=c_nkn$ in \eqref{lge2}, where $\lim_{n\to\infty}c_n=c$,   and letting $n\to\infty$ and recalling
\eqref{keylimit}, which was derived for $s\ge1$, but which holds just as well with $s=0$, it follows  that
\begin{equation}\label{l2ton}
\begin{aligned}
&\lim_{n\to\infty}P_{n,k}(J_{n,k}^{\mathcal{S}}=l|W_{\mathcal{S}(n,k;M_n)})=
\frac1{\lim_{n\to\infty}P_{n,k}(W_{\mathcal{S}(n,k;M_n)})}  \binom k{l-1}(\frac c{1-c})^{l-1}(1-c)^{k}=\\
&\frac1{\lim_{n\to\infty}P_{n,k}(W_{\mathcal{S}(n,k;M_n)})}  \binom k{l-1}c^{l-1}(1-c)^{k-l+1},\  l\in\{2,\cdots, k\},
\end{aligned}
\end{equation}
which gives \eqref{pi} for $l\in\{2,\cdots, k\}$.
Since
\begin{equation}\label{binomial0k}
\sum_{l=2}^k \binom k{l-1}c^{l-1}(1-c)^{k-l+1}=1-c^k-(1-c)^k,
\end{equation}
it follows from \eqref{l2ton} and \eqref{binomial0k} that
\eqref{pi} also holds for $l=1$.
\hfill $\square$

\section{Proof of Theorem \ref{1tokn}}\label{quick2}
\it\noindent  Proof of Theorem \ref{1tokn}.\rm\
Recall that $A^{(n)}_{M,j,l}$ denotes the event that among the first $M$ items,
 the number $j$ has occurred $l$ times, and no number greater than $j$ has occurred at all,
 where $1\le M\le kn-1,\ j\in[n]$ and  $l\in[k]$. We first consider the case that
the event $A^{(n)}_{M,j,l}$  occurs, where $j\in[n-1]$,
and then afterwards, we consider the case that the event  $A^{(n)}_{M,n,l}$ occurs.

Let $j\in[n-1]$. If $A^{(n)}_{M,j,l}$  occurs, then among the last $kn-M$ items, there are
$(j-1)k-M+l$ items of rank strictly less than $j$, there are $k(n-j)-l$ items of ranks between $j$ and $n-1$, and there are $k$ items of rank $n$.

On the event
$W_{\mathcal{S}(n,k;M)}\cap A^{(n)}_{M,j,l}$, the intersection of the event $A^{(n)}_{M,j,l}$ and the event  that the strategy $\mathcal{S}(n,k;M)$ successfully selects an
item of highest rank $n$, let $\tau(A^{(n)}_{M,j,l})$ denote the number of additional items after the first $M$ items  up until an item of rank $n$ is selected.
Then by definition,
\begin{equation}\label{X+c}
X_{n,k}^{\mathcal{S}}=M+\tau(A^{(n)}_{M,j,l})\ \text{on}\ W_{\mathcal{S}(n,k;M)}\cap A^{(n)}_{M,j,l}.
\end{equation}
From the previous paragraph and from the definition of the strategy $\mathcal{S}(n,k;M)$,  the distribution of $\tau(A^{(n)}_{M,j,l})$ under
 $P_{n,k}(\cdot\ |W_{\mathcal{S}(n,k;M)}\cap A^{(n)}_{M,j,l})$ is the same as the distribution
of $\tau^{(N,A)}$ in Proposition \ref{combibox} with $N=(j-1)k-M+l$ and $A=k$.
For convenience, we copy formula \eqref{condprobj<n} here:
\begin{equation}\label{WcondA}
P_{n,k}(W_{\mathcal{S}(n,k;M)}|A^{(n)}_{M,j,l})=\frac k{k(n-j+1)-l},\ j\in[n-1].
\end{equation}

Now consider the case that the event $A^{(n)}_{M,n,l}$  occurs.
 Then among the last $kn-M$ items, there are
$(n-1)k-M+l$ items of rank strictly less than $n$, and there are $k-l$ items of rank $n$.

On the event
$W_{\mathcal{S}(n,k;M)}\cap A^{(n)}_{M,n,l}$, the intersection of the event $A^{(n)}_{M,n,l}$ and the event  that the strategy $\mathcal{S}(n,k;M)$ successfully selects an
item of highest rank $n$, let $\tau(A^{(n)}_{M,n,l})$ denote the number of additional items after the first $M$ items  up until an item of rank $n$ is selected.
Then by definition,
\begin{equation}\label{X+cagain}
X_{n,k}^{\mathcal{S}}=M+\tau(A^{(n)}_{M,n,l})\ \text{on}\ W_{\mathcal{S}(n,k;M)}\cap A^{(n)}_{M,n,l}.
\end{equation}
From the previous paragraph and from the definition of the strategy $\mathcal{S}(n,k;M)$,  the distribution of $\tau(A^{(n)}_{M,n,l})$ under
 $P_{n,k}(\cdot\ |W_{\mathcal{S}(n,k;M)}\cap A^{(n)}_{M,n,l})$ is the same as the distribution
of $\tau^{(N,A)}$ in Proposition \ref{combibox} with $N=(n-1)k-M+l$ and $A=k-l$.
For convenience, we copy formula \eqref{condprob=n} here:
\begin{equation}\label{WcondAn}
P_{n,k}(W_{\mathcal{S}(n,k;M)}|A^{(n)}_{M,n,l})=\begin{cases}1,\ \text{if}\ l\in[k-1];\\ 0,\ \text{if}\ l=k.\end{cases}
\end{equation}

Let $\{M_n\}$ satisfy $M_n= c_nkn$ with $\lim_{n\to\infty}c_n=c\in(0,1)$. Let $H(n,k;M_n)$ denote  the expected number of additional items
observed after the first $M_n$ items, conditioned on successfully selecting an item of rank $n$ while using strategy $\mathcal{S}(n,k;M_n)$.
For ease of notation, let
$$
\begin{aligned}
&\beta_{n,k,M_n}:=P_{n,k}(W_{\mathcal{S}(n,k;M_n)});\\
&\mathcal{WA}^{(n,k)}_{M_n,j,l}:=W_{\mathcal{S}(n,k;M_n)}\cap A^{(n)}_{M_n,j,l}.
\end{aligned}
$$
From the above analysis, it follows that
$$
\begin{aligned}
&H(n,k;M_n)=\\
&\frac1{\beta_{n,k,M_n}}\sum_{j=1}^{n-1}\sum_{l=1}^kP(\mathcal{WA}^{(n,k)}_{M_n,j,l})E_{n,k}(\tau(A^{(n)}_{M_n,j,l})|\mathcal{WA}^{(n,k)}_{M_n,j,l})+\\
&\frac1{\beta_{n,k,M_n}}\sum_{l=1}^kP(\mathcal{WA}^{(n,k)}_{M_n,n,l})E_{n,k}(\tau(A^{(n)}_{M_n,n,l})|\mathcal{WA}^{(n,k)}_{M_n,n,l})=\\
&\frac1{\beta_{n,k,M_n}}\sum_{j=1}^{n-1}\sum_{l=1}^kP(A^{(n)}_{M_n,j,l})\frac k{k(n-j+1)-l}E_{n,k}(\tau(A^{(n)}_{M_n,j,l})|\mathcal{WA}^{(n,k)}_{M_n,j,l})+\\
&\frac1{\beta_{n,k,M_n}}\sum_{l=1}^{k-1}P(A^{(n)}_{M_n,n,l})E_{n,k}(\tau(A^{(n)}_{M_n,n,l})|\mathcal{WA}^{(n,k)}_{M_n,n,l}),
\end{aligned}
$$
where \eqref{WcondA} and \eqref{WcondAn} were used for the second equality.
Substituting from \eqref{A} in the right hand side of the  above equation, and substituting $M_n=c_nkn$, we obtain after division by $nk$,
\begin{equation}\label{Hn}
\begin{aligned}
&\frac1{nk}H(n,k;M_n)=
\frac1{\beta_{n,k,c_nkn}}\times\\
&\sum_{j=1}^{n-1}\sum_{l=1}^k\frac{\binom {c_nkn}l(k)_l(k(j-1))_{c_nkn-l}}{(kn)_{c_nkn}}\frac k{k(n-j+1)-l}
\frac{E_{n,k}(\tau(A^{(n)}_{c_nkn,j,l})|\mathcal{WA}^{(n,k)}_{M_n,j,l})}{nk}+\\
&\frac1{\beta_{n,k,c_nkn}}\sum_{l=1}^{k-1}\frac{\binom {c_nkn}l(k)_l(k(n-1))_{c_nkn-l}}{(kn)_{c_nkn}}\frac{E_{n,k}(\tau(A^{(n)}_{c_nkn,n,l})|\mathcal{WA}^{(n,k)}_{M_n,n,l})}{nk}.
\end{aligned}
\end{equation}

Consider first  the second sum on the right hand side of \eqref{Hn}.
From earlier in the proof, under $P_{n,k}(\cdot\ |\mathcal{WA}^{(n,k)}_{M_n,n,l})$, the random variable
$\tau(A^{(n)}_{c_nkn,n,l})$ has the distribution  of $\tau^{N,A}$  from Proposition \ref{combibox} with $N=(n-1)k-c_nkn+l$ and $A=k-l$.
Thus, by Proposition \ref{combibox},
\begin{equation}\label{taulimn}
\lim_{n\to\infty}\frac{E_{n,k}(\tau(A^{(n)}_{c_nkn,n,l})|\mathcal{WA}^{(n,k)}_{M_n,n,l})}{nk}=\frac{1-c}{k-l+1}.
\end{equation}
Using \eqref{taulimn} and
\eqref{finalfirsttermformula}, it follows that the second term on the right hand side of \eqref{Hn} satisfies
\begin{equation}\label{keyagain2}
\begin{aligned}
&\lim_{n\to\infty}\frac1{\beta_{n,k,c_nkn}}\sum_{l=1}^{k-1}\frac{\binom {c_nkn}l(k)_l(k(n-1))_{c_nkn-l}}{(kn)_{c_nkn}}\frac{E_{n,k}(\tau(A^{(n)}_{c_nkn,n,l})|\mathcal{WA}^{(n,k)}_{M_n,n,l})}{nk}=\\
&\frac1{\lim_{n\to\infty}P_{n,k}(W_{\mathcal{S}(n,k;M_n)})}(1-c)^k\sum_{l=1}^{k-1}\binom kl(\frac c{1-c})^l\thinspace\frac{1-c}{k-l+1}=\\
&\frac1{\lim_{n\to\infty}P_{n,k}(W_{\mathcal{S}(n,k;M_n)})}(1-c)\sum_{i=2}^k\binom k{i-1}c^{i-1}(1-c)^{k-i+1}\frac1{k+2-i}=\\
&(1-c)\sum_{i=2}^k\frac{p_c^{(k)}(i)}{k+2-i},
\end{aligned}
\end{equation}
where the last step follows from \eqref{pi}.

Consider now the first sum on the right hand side of \eqref{Hn}.
For $N\in\mathbb{N}$, we break the outer sum $\sum_{j=1}^{n-1}$  into two pieces--the sum from 1 to $n-N$, and the sum from $n-N+1$ to $n-1$.
Since $\frac{E_{n,k}(\tau(A^{(n)}_{c_nkn,j,l})|\mathcal{WA}^{(n,k)}_{M_n,j,l})}{nk}$ is bounded, it follows from
\eqref{Nnegligible} that the part of the sum from 1 to $n-N$ satisfies
\begin{equation}\label{Nnegligibleagain}
\begin{aligned}
&\lim_{N\to\infty}\lim_{n\to\infty}\\
&\sum_{j=1}^{n-N}\sum_{l=1}^k\frac{\binom {c_nkn}l(k)_l(k(j-1))_{c_nkn-l}}{(kn)_{c_nkn}}\frac k{k(n-j+1)-l}\frac{E_{n,k}(\tau(A^{(n)}_{c_nkn,j,l})|\mathcal{WA}^{(n,k)}_{M_n,j,l})}{nk}=0.
\end{aligned}
\end{equation}
Letting $s=n-j$, the part of the sum from $n-N+1$ to $n-1$ becomes
\begin{equation}\label{nonnegl}
\begin{aligned}
&\sum_{j=n-N+1}^{n-1}\sum_{l=1}^k\frac{\binom {c_nkn}l(k)_l(k(j-1))_{c_nkn-l}}{(kn)_{c_nkn}}\frac k{k(n-j+1)-l}\frac{E_{n,k}(\tau(A^{(n)}_{c_nkn,j,l})|\mathcal{WA}^{(n,k)}_{M_n,j,l})}{nk}=\\
&\sum_{s=1}^{N-1}\sum_{l=1}^k\frac{\binom {c_nkn}l(k)_l(k(n-s-1))_{c_nkn-l}}{(kn)_{c_nkn}}\frac k{k(s+1)-l}\frac{E_{n,k}(\tau(A^{(n)}_{c_nkn,n-s,l})|\mathcal{WA}^{(n,k)}_{M_n,n-s,l})}{nk}.
\end{aligned}
\end{equation}
Recall from earlier in the proof that under $P_{n,k}(\cdot\ |\mathcal{WA}^{(n,k)}_{c_nkn,n-s,l})$, the random variable
$\tau(A^{(n)}_{c_nkn,n-s,l})$ has the distribution  of $\tau^{N,A}$  from Proposition \ref{combibox} with $N=(n-s-1)k-c_nkn+l$ and $A=k$.
Thus, by Proposition \ref{combibox},
\begin{equation}\label{taulim}
\lim_{n\to\infty}\frac{E_{n,k}(\tau(A^{(n)}_{c_nkn,n-s,l})|\mathcal{WA}^{(n,k)}_{M_n,n-s,l})}{nk}=\frac{1-c}{k+1}.
\end{equation}
Using  \eqref{nonnegl}, \eqref{taulim} and
\eqref{keylimit} (see also \eqref{keytermlimit}), it follows that the first term on the right hand side of \eqref{Hn} satisfies
\begin{equation}\label{keyagain1}
\begin{aligned}
&\lim_{n\to\infty}\frac1{\beta_{n,k,c_nkn}}\times\\
&\sum_{j=n-N+1}^{n-1}\sum_{l=1}^k\frac{\binom {c_nkn}l(k)_l(k(j-1))_{c_nkn-l}}{(kn)_{c_nkn}}\frac k{k(n-j+1)-l}\frac{E_{n,k}(\tau(A^{(n)}_{c_nkn,j,l})|\mathcal{WA}^{(n,k)}_{M_n,j,l})}{nk}=\\
&\frac{1-c}{k+1}\Big(\frac1{\lim_{n\to\infty}P_{n,k}(W_{\mathcal{S}(n,k;M_n)})}\sum_{l=1}^k\binom kl(\frac c{1-c})^l\sum_{s=1}^{N-1}(1-c)^{k(s+1)}\frac k{k(s+1)-l}\Big).
\end{aligned}
\end{equation}

A perusal of   the current proof,  or alternatively, a perusal of the proofs of Proposition \ref{prop1} and  Theorem \ref{thm1},   will reveal that
the expression in the large parentheses on the right hand side of \eqref{keyagain1} must
be equal to $p_c^{(k)}(1)$. (In the statement and proof of Theorem \ref{1tok}, we calculated $p_c^{(k)}(1)$ more simply as
$1-\sum_{i=2}^kp_c^{(k)}(i)$.)
Substituting this in \eqref{keyagain1}, and using \eqref{Nnegligibleagain}, it follows that
\begin{equation}\label{firstsum}
\begin{aligned}
&\lim_{n\to\infty}\sum_{j=1}^{n-1}\sum_{l=1}^k\frac{\binom {c_nkn}l(k)_l(k(j-1))_{c_nkn-l}}{(kn)_{c_nkn}}\frac k{k(n-j+1)-l}
\frac{E_{n,k}(\tau(A^{(n)}_{c_nkn,j,l})|\mathcal{WA}^{(n,k)}_{M_n,j,l})}{nk}=\\
&(1-c)\frac{p_c^{(k)}(1)}{k+1}.
\end{aligned}
\end{equation}
The theorem in the form provided by the first equality in \eqref{gamma}  follows from  \eqref{X+c}, \eqref{X+cagain}, the definition of $H(n,k;M_n)$ in the paragraph following
\eqref{WcondAn},   \eqref{Hn}, \eqref{keyagain2} and  \eqref{firstsum}. The second equality in \eqref{gamma} follows easily from the definition of
$\{p_c^{(k)}(i)\}_{i=1}^n$ in \eqref{pi}.
\hfill $\square$

\section{Proofs of Proposition \ref{prop+} and Theorem \ref{thm+}}\label{sec4}
\noindent \it Proof of Proposition \ref{prop+}.\rm\  We follow the proof of Proposition \ref{prop1} up through \eqref{A}.
Consider now the case that $j\in[n-1]$. If one employs the strategy $\mathcal{S}^+(n,k;M)$, and the event $A^{(n)}_{M,j,l}$ occurs, then after the $M$th items arrives, one will select the first item strictly larger than $j$. Among the items arriving after
the $M$th item arrives, there are $k(n-j)$ items that are larger than $j$, of which $k$ of them have the highest rank $n$. Since items arrive in uniformly random order, the conditional probability of
selecting an item of highest rank, conditioned on  $A^{(n)}_{M,j,l}$,  is given by
\begin{equation}\label{condprobj<n+}
P_{n,k}(W_{\mathcal{S}^+(n,k;M)}|A^{(n)}_{M,j,l})=\frac 1{(n-j)},\ j\in[n-1].
\end{equation}

 On the other hand, consider now the case  $j=n$. If one  employs the strategy $\mathcal{S}^+(n,k;M)$, and the event $A^{(n)}_{M,n,l}$ occurs, then after the $M$th item arrives,  one is required to
 select the first item that is greater than  $n$, but of course  no such item exists; therefore,
 \begin{equation}\label{condprob=n+}
  P_{n,k}(W_{\mathcal{S}^+(n,k;M)}|A^{(n)}_{M,n,l})=0.
\end{equation}

 Since for each $M$,  the collection of events $\{A^{(n)}_{M,j,l}: j\in[n], l\in [k]\thinspace\}$ are disjoint, and the probability of their union is 1, it follows from \eqref{A}, \eqref{condprobj<n+} and \eqref{condprob=n+} that
 \begin{equation*}
 P_{n,k}(W_{\mathcal{S}^+(n,k;M)})=
\sum_{j=1}^{n-1}\sum_{l=1}^k \frac{\binom Ml(k)_l(k(j-1))_{M-l}}{(kn)_M}\frac 1{(n-j)}.
 \end{equation*}
 \hfill $\square$

\noindent \it Proof of Theorem \ref{thm+}.\rm\ The right hand side of \eqref{formula1+}, with $M=M_n$, is equal to the expression whose limit is taken on the  left hand side of
\eqref{finalsecondtermformula}, but with the fraction $\frac1{k(n-j+1)-l}$ occurring there  replaced by $\frac1{k(n-j)}$. Consequently,
\begin{equation}\label{limit+}
\begin{aligned}
&\lim_{n\to\infty}P_{n,k}(W_{\mathcal{S}^+(n,k;M_n)})=\Big(k\sum_{l=1}^k\binom kl(\frac c{1-c})^l\Big)\Big(\sum_{s=1}^\infty(1-c)^{k(s+1)}\frac 1{ks}\Big),\\
& \text{for}\  M_n\sim cn, \ c\in(0,1),
\end{aligned}
\end{equation}
where the right  hand side above is the right hand side of \eqref{finalsecondtermformula}, but  with the fraction $\frac1{k(s+1)-l}$ occurring there replaced by
$\frac1{ks}$.  Note that the infinite series occurring on the right hand side of \eqref{limit+} is equal to $G_{k,k}(1-c)$, which was defined in \eqref{G}. Thus, it follows from
\eqref{limit+} and \eqref{Gformula} that
\begin{equation}\label{limit+final}
\begin{aligned}
&\lim_{n\to\infty}P_{n,k}(W_{\mathcal{S}^+(n,k;M_n)})=\Big(k\sum_{l=1}^k\binom kl(\frac c{1-c})^l\Big)\Big((-\frac{(1-c)^k}k\log(1-(1-c)^k)\Big)=\\
&\Big((\frac c{1-c}+1)^k-1\Big)\Big(-(1-c)^k\log(1-(1-c)^k)\Big)=\\
&-\big(1-(1-c)^k\big)\log(1-(1-c)^k).
\end{aligned}
\end{equation}
\hfill$\square$

\end{document}